\documentclass[11pt]{amsart}

\usepackage{amsmath,amssymb,amsthm,amsfonts,verbatim}
\usepackage{enumerate}
\usepackage{microtype}
\usepackage{epsfig}
\usepackage[font=small,format=plain,labelfont=bf,up,textfont=it,up]{caption}
\usepackage{pinlabel}
\usepackage[all,2cell]{xy}
\usepackage{mathtools}
\usepackage{url}

\usepackage[top=0.9in,left=1in,right=1in,bottom=1in]{geometry}

\title{The rational cohomology of the mapping class group\\ vanishes
 in its virtual cohomological dimension}

\author{Thomas Church, Benson Farb and Andrew Putman }
\thanks{The authors gratefully acknowledge  support from the National Science Foundation.}

\theoremstyle{plain}
\newtheorem{theorem}{Theorem}
\theoremstyle{definition}

\newcommand{\nc}{\newcommand}
\nc{\dmo}{\DeclareMathOperator}

\nc{\Q}{\mathbb{Q}}
\nc{\R}{\mathbb{R}}
\nc{\Z}{\mathbb{Z}}
\dmo{\Mod}{Mod}
\dmo{\SL}{SL}
\nc{\M}{\mathcal{M}}
\nc{\T}{\mathcal{T}}
\nc{\C}{\mathcal{C}}
\nc{\U}{\mathcal{U}}
\nc{\A}{\mathcal{A}}
\dmo{\thick}{thick}
\dmo{\vcd}{vcd}
\dmo{\cd}{cd}
\dmo{\St}{St}

\nc{\Tthick}{\T_g^{\thick}}
\nc{\Mthick}{\M_g^{\thick}}

\renewcommand{\epsilon}{\varepsilon}
\nc{\coloneq}{\mathrel{\mathop:}\mkern-1.2mu=}
\nc{\margin}[1]{\marginpar{\scriptsize #1}}
\nc{\para}[1]{\medskip\noindent\textbf{#1.}}

\newcommand\Figure[3]{
\begin{figure}[t]
\centering
\centerline{\psfig{file=#2,scale=60}}
\caption{#3}
\label{#1}
\end{figure}}

\begin{document}

\begin{abstract}
Let $\Mod_g$ be the mapping class group of a genus $g\geq 2$ surface.  The
group $\Mod_g$ has virtual cohomological dimension $4g-5$.  In this note we use a theorem of
Broaddus and the combinatorics of chord diagrams to prove that $H^{4g-5}(\Mod_g;\Q)=0$.
\end{abstract}
\maketitle

\section{Introduction}
Let $\Mod_g$ be the mapping class group of a closed, oriented, genus $g\geq 2$ surface, and let ${\mathcal M}_g$ 
be the moduli space of genus $g$ Riemann surfaces.  It is well-known 
that for each $i\geq 0$,
\[H^i(\Mod_g;\Q) \cong H^i({\mathcal M}_g;\Q).\]
It is a fundamental open problem to determine the maximal $i$ for which these vector 
spaces are nonzero.  Harer \cite{Ha} proved that the \emph{virtual cohomological dimension} 
$\vcd(\Mod_g)$ equals $4g-5$.  More precisely, he proved that $H^{4g-5}(\Mod_g;\St_g \otimes\, \Q) \neq 0$
for a certain $\Mod_g$-module $\St_g$ (see below for details) and 
that $H^{i}(\Mod_g;V \otimes \Q)=0$ for all $i > 4g-5$ and all $\Mod_g$-modules $V$.
 Thus the first step of the problem above is to determine whether $H^{4g-5}(\Mod_g;\Q)\neq 0$.  The purpose of this note is to answer this question.

Let $\Mod_{g,\ast}$ (resp.\ $\Mod_{g,1}$) denote the mapping class group of the genus $g$ surface with one marked point (resp.\ one boundary component).
\begin{theorem}\label{theorem:main}
For any $g\geq 2$,
\[H^{4g-5}(\Mod_g;\Q)=H^{4g-5}({\mathcal M}_g;\Q)=0.
\]
Further, the rational cohomology of $\Mod_{g,\ast}$ (resp.\ the integral cohomology of $\Mod_{g,1}$) vanishes in its virtual cohomological dimension. 
\end{theorem}

This theorem was announced some years ago by Harer, but he has informed us that his proof will not appear. We recently learned that Morita--Sakasai--Suzuki \cite{MSS} have independently found a proof of Theorem~\ref{theorem:main} using a completely different method.  They apply a theorem of Kontsevich on graph homology to their computation of a generating set for a certain symplectic Lie algebra.   Our proof combines some results about
the combinatorics of chord diagrams with the work of Broaddus \cite{Br}
on the Steinberg module of $\Mod_g$. We thank Allen Hatcher and Takuya Sakasai for their comments on an earlier version of this paper, and John Harer for informing us about the paper \cite{MSS} and his own work.

Theorem~\ref{theorem:main} is consistent with the well-studied analogy between mapping class groups and arithmetic groups. For example, Theorem~1.3 of Lee--Szczarba \cite{LS} states that the rational cohomology of $\SL(n,\Z)$ vanishes in its cohomological dimension.

\section{Background}
We begin by briefly summarizing previous results that make our computation possible; for details see Broaddus \cite{Br}.

\para{Teichm\"{u}ller space and its boundary} Let $S_g$ be a connected, closed orientable surface of genus $g\geq 2$.
 Let $\C_g$ be the \emph{curve complex} of $S_g$ defined by Harvey \cite{Harv}, i.e.\ the flag complex whose $k$-simplices are the $(k+1)$-tuples of distinct free homotopy classes of simple closed curves in $S_g$ that can be realized disjointly.  Harer \cite{Ha} proved that $\C_g$ is homotopy equivalent to a wedge of spheres $\bigvee_{i=1}^\infty S^{2g-2}$.

There exists a constant $\delta>0$ such that any two closed geodesics on a hyperbolic surface of length $\leq \delta$ are disjoint (the \emph{Margulis constant} for hyperbolic surfaces).  Let $\Tthick$ be the Teichm\"uller space of marked hyperbolic surfaces diffeomorphic to $S_g$ having no closed geodesic of length $<\delta$. It is known that $\Tthick$ is a $(6g-6)$-dimensional manifold with corners. Ivanov \cite{Iv} proved that $\mathcal{T}_g^{\thick}$ is contractible and that its boundary $\partial \Tthick$ is homotopy equivalent to $\C_g$. Briefly, for each simplex $\sigma$ of $\C_g$, let $\T_{\sigma}$ be the subset of 
$\partial \Tthick$ consisting of surfaces where each curve in $\sigma$ has length $\delta$.
Each $T_{\sigma}$ is contractible, and $\T_{\sigma} \cap \T_{\sigma'} = \emptyset$ unless $\sigma \cup \sigma'$ is a simplex of $\C_g$, in which case $T_{\sigma} \cap \T_{\sigma'} = \T_{\sigma \cup \sigma'}$.

\para{Duality in the mapping class group} The mapping class group $\Mod_g$ acts properly discontinuously on $\Tthick$ with finite stabilizers.  Defining 
$\Mthick = \Tthick / \Mod_g$, it follows that 
$H^*(\Mod_g;\Q)\cong H^*(\Mthick;\Q)$. Mumford's compactness criterion states that $\Mthick$ is compact. 
Combining this with the previous two paragraphs, the work of Bieri--Eckmann \cite[Theorem~6.2]{BE} shows that 
$\vcd(\Mod_g)=4g-5$ and that
\begin{equation}
\label{eq:duality}
H^{4g-5}(\Mod_g;\Q)\cong H_0(\Mod_g;H_{2g-2}(\C_g;\Q)).
\end{equation}

In fact, we can say more.  Let $\St_g$ denote the \emph{Steinberg module}, i.e.\ the $\Mod_g$-module
$H_{2g-2}(\C_g;\Z)$.  Then $\St_g\otimes\, \Q$ is the rational \emph{dualizing module} for $\Mod_g$, meaning that 
\[H^{4g-5-k}(\Mod_g;M\otimes\Q)\cong H_k(\Mod_g;M\otimes \St_g\otimes\, \Q)\] 
for any $k$ and any $M$. Moreover $\St_g$ is also the dualizing module for $\Mod_{g,*}$ and $\Mod_{g,1}$,
which act on $\St_g$ via the natural surjections $\Mod_{g,*}\to \Mod_g$ and $\Mod_{g,1}\to \Mod_g$ \cite{Ha}. This implies that for $\nu=\vcd(\Mod_{g,*})=4g-3$ we have $H^{\nu-k}(\Mod_{g,*};M\otimes\Q)\cong H_k(\Mod_{g,*};M\otimes \St_g\otimes\, \Q)$. For $\Mod_{g,1}$ we obtain a similar result with $\nu=\cd(\Mod_{g,1})=4g-2$, except that since $\Mod_{g,1}$ is torsion-free the result holds integrally: $H^{\nu-k}(\Mod_{g,1};M)\cong H_k(\Mod_{g,1};M\otimes \St_g)$.

\para{An alternate model for $\St_g$} Fix a finite-volume hyperbolic metric on $S_g-\{\ast\}$. Another model for $\St_g$ comes from the \emph{arc complex} $\A_g$, the flag complex whose $k$-simplices are the disjoint $(k+1)$-tuples of simple geodesics on $S_g-\{\ast\}$ beginning and ending at the cusp $\ast$. Let $\A_g^\infty$ be the subcomplex consisting of collections of geodesics $\gamma_1,\ldots,\gamma_{k+1}$ for which $S-\bigcup \gamma_i$ has some non-contractible component. Harer proved that $\A_g^\infty$ is homotopy equivalent to $\C_g$ \cite{Ha}, and that $\A_g$ is contractible \cite{Ha2} (see also \cite{Hat}). Thus
\[\St_g=H_{2g-2}(\C_g)\simeq H_{2g-2}(\A_g^\infty)\simeq H_{2g-1}(\A_g/\A_g^\infty).\]

\para{Chord diagrams} By examining how the geodesics are arranged in a neighborhood of $\ast$, an $(n-1)$-simplex of $\A_g$ can be encoded by a $n$-chord diagram; see \cite[\S4.1]{Br}. An \emph{ordered $n$-chord diagram} is an ordered sequence $U = (u_1,\ldots,u_n)$, where $u_i$ is an unordered pair of 
distinct points on $S^1$ (a \emph{chord}) and
$u_i \cap u_j = \emptyset$ if $i \neq j$.
 We will visually depict $U$ by drawing arcs connecting
the points in each $u_i$ (see Figure \ref{figure:1} for examples).
Two ordered chord diagrams are identified if they differ by an orientation-preserving
homeomorphism of the circle.  

\para{Filling systems} An unlabeled \emph{$k$-filling system} of genus $g$ is a $(2g+k)$-chord diagram satisfying the conditions described in \cite[\S4.1]{Br}: no chord should be parallel to another chord or to the boundary circle, and the chords should determine exactly $k+1$ boundary cycles. These conditions, which guarantee that these chords define a simplex of $\A_g-\A_g^\infty$, have the following simple combinatorial formulation. Given $U = (u_1,\ldots,u_n)$, consider two permutations of the $2n$ points $u_1\cup\cdots\cup u_n$: let $\omega$ be the $2n$-cycle which takes each point to the point immediately adjacent in the clockwise direction, while $\tau$ exchanges the two points of each chord $u_i$ and thus is a product of $n$ transpositions. Then a $(2g+k)$-chord diagram is a $k$-filling system of genus $g$ if $\tau\circ \omega$ has $k+1$ orbits, none of which have length 1 or 2. Finally, let $t_i$ be the straight line in $D^2$ connecting the two points of $u_i$. Then we say that $U$ is \emph{disconnected} if the set $t_1\cup \cdots\cup t_n\subset D^2$ is not connected.

\para{The chord diagram chain complex}
Fix a genus $g$, and set $n=2g+k$.  Let $\U_k$ be the free abelian group spanned by ordered $k$-filling systems of genus
$g$ modulo the following relation.  For $\sigma\in S_n$ 
and $U = (u_1,\ldots,u_n)$, define $\sigma\cdot U= (u_{\sigma(1)},\ldots,u_{\sigma(n)})$. We impose the relation $\sigma\cdot U = (-1)^{\sigma} U$.
The differential $\partial : \U_k \rightarrow \U_{k-1}$ is defined as follows.  Consider
an ordered $k$-filling system $U = (u_1,\ldots,u_n)$ of genus $g$.  For $1 \leq i \leq n$, let
$\partial_i U$ equal $(u_1,\ldots,\widehat{u_i},\ldots,u_n)$ if this is an ordered
$(k-1)$-filling system of genus $g$; otherwise, let $\partial_i U = 0$.
Then \[\partial(U) = \sum_{i=1}^n (-1)^{i-1} \partial_i U.\]

\para{Broaddus's results}
We will need the following theorem of Broaddus \cite{Br}.  Recall that if $\Gamma$ is a group
and $M$ is a $\Gamma$-module, then the \emph{module of coinvariants}, denoted $M_{\Gamma}$, is
the quotient $M / \langle g \cdot m - m\,|\,g \in \Gamma, m \in M\rangle$. Let $X$ be the $0$-filling system of genus $g$ depicted in Figure \ref{figure:1}a.

\begin{theorem}[{Broaddus \cite{Br}}]
\label{theorem:broaddus}For each $g\geq 0$, the following hold.
\begin{enumerate}[(i)]
\item $(\St_g)_{\Mod_g} \cong \U_0 / \partial(\U_1)$.
\item The abelian group $\U_0 / \partial(\U_1)$
is spanned by the image $[X]\in \U_0/\partial(\U_1)$ of $X\in \U_0$.
\item If $v$ is a disconnected $0$-filling system of genus $g$, then the image
of $v$ in $\U_0 / \partial(\U_1)$ is $0$.
\end{enumerate}
\end{theorem}

\noindent
For part (i) of Theorem \ref{theorem:broaddus}, see \cite[Proposition 3.3]{Br} together with the
remark preceding \cite[Example 4.1]{Br}; for part (ii), see
\cite[Theorem 4.2]{Br}; and for part (iii), see
\cite[Proposition 4.5]{Br}.

\section{Proof of Theorem~\ref{theorem:main}}
For any group $\Gamma$ and any $\Gamma$-module $M$, recall that $H_0(\Gamma;M) = M_{\Gamma}$.
Since the actions of $\Mod_{g,\ast}$ and $\Mod_{g,1}$ on $\St_g$ factor through $\Mod_g$, to 
prove Theorem~\ref{theorem:main} it suffices by \eqref{eq:duality} 
to show that $(\St_g)_{\Mod_g} = 0$.
By Theorem \ref{theorem:broaddus}(i), this is equivalent to showing that
$\U_0 / \partial(\U_1) = 0$.  

For $v \in \U_0$, let $[v]$ denote the associated element of $\U_0 / \partial(\U_1)$.
Let $X = (x_1,\ldots,x_{2g})$ be the $0$-filling system depicted in Figure \ref{figure:1}(a).
By Theorem \ref{theorem:broaddus}(ii), it is enough to show that $[X]=0$.
Let $Y = (x_1,\ldots,x_{2g},y)$ be the $1$-filling system depicted in Figure \ref{figure:1}(b). 
Observe that
\[\partial_1 Y = (x_2,\ldots,x_{2g},y) = (x_1,\ldots,x_{2g}) = X,\]
where the second equality holds since the indicated chord diagrams differ by
an orientation preserving homeomorphism of $S^1$.  Similarly, $\partial_{2g+1}Y = X$.
Also, $\partial_2 Y = 0$ (resp.\ $\partial_{2g}Y = 0$) by definition, since the chord $x_1$ (resp.\ $x_{2g+1}$) becomes parallel to the boundary.  We thus have
\[\partial(Y) = 2X + \sum_{i=3}^{2g-1} (-1)^{i-1} \partial_i Y.\]
For $3 \leq i \leq 2g-1$, the $0$-filling system $\partial_i Y$ is disconnected, so Theorem \ref{theorem:broaddus}(iii) implies that $[\partial_i Y]=0$.  We conclude that
$2[X]=0$.

\Figure{figure:1}{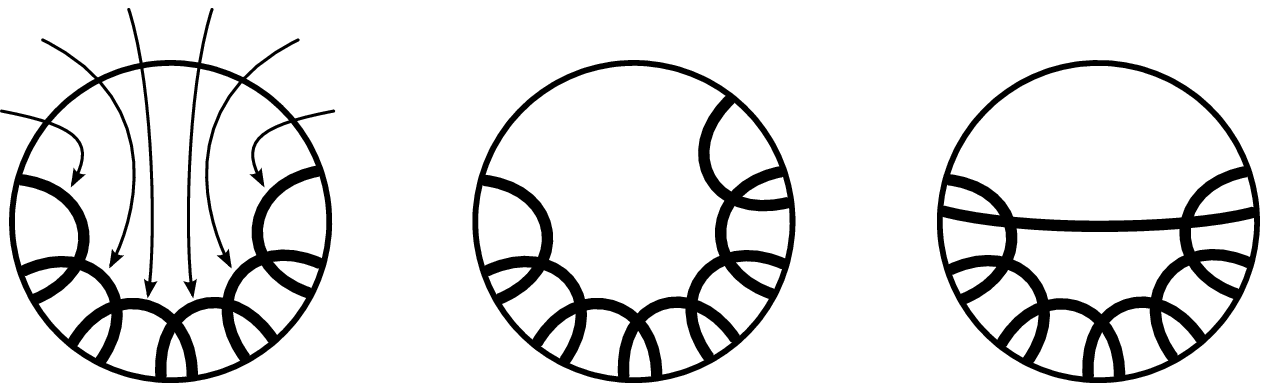}{(a) The oriented $0$-filling system $X = (x_1,\ldots,x_{2g})$.  For concreteness, we depict
it for $g = 3$.  In general, $X$ has $2g$ chords arranged in the same pattern as the chords shown.\\
 (b) The $1$-filling system $Y = (x_1,\ldots,x_{2g},y)$.  The chord $y$ intersects the chord $x_{2g}$.\\
(c) The $1$-filling system $Z = (z,x_1,\ldots,x_{2g})$.  The chord $z$ intersects both $x_1$ and $x_{2g}$.}

Now consider the $1$-filling system $Z = (z,x_1,\ldots,x_{2g})$ depicted in Figure \ref{figure:1}(c).
Removing any chord from Figure \ref{figure:1}(c) yields Figure \ref{figure:1}(a) up to rotation, 
so $\partial_i Z=\pm X$ for each $i$.  In fact, it is clear that $\partial_1 Z=X$, that $\partial_2 Z = -X$, that
$\partial_3 Z = X$, and so on, with $\partial_i Z = (-1)^{i-1} X$.  This shows that
\[\partial (Z) = X + X + \cdots X = (2g+1)X,\]
so $(2g+1)[X] = 0$.

Summing up, we have shown that $2[X] = (2g+1)[X] = 0$.  This implies that $[X]=0$, as desired.

\small

\noindent
E-mail: tchurch@math.uchicago.edu, farb@math.uchicago.edu, andyp@rice.edu


\begin{thebibliography}{ABCDEF}
\small

\bibitem[BE]{BE} R. Bieri and B. Eckmann, Groups with homological duality generalizing Poincar\'{e} duality, \emph{Invent. Math.} 20 (1973), 103--124.

\bibitem[Br]{Br} N. Broaddus, Homology of the curve complex and the Steinberg module of the mapping class group, preprint (2007), arXiv:0711.0011v2.



\bibitem[Ha]{Ha} J. Harer, The virtual cohomological dimension of the mapping class group of an orientable surface, \emph{Invent. Math.}  84  (1986),  no. 1, 157--176.

\bibitem[Ha2]{Ha2} J. Harer, Stability of the homology of the mapping class groups of orientable surfaces, \emph{Ann. of Math.}  121  (1985),  no. 2, 215--249.

\bibitem[Harv]{Harv} W. J. Harvey, Boundary structure of the modular group, in
 \emph{Riemann surfaces and related topics: Proceedings of the 1978 Stony
 Brook Conference}, 245--251, Ann. of Math. Stud., 97, Princeton Univ. Press, Princeton, N.J.,  1981.

\bibitem[Hat]{Hat} A. Hatcher, On triangulations of surfaces, \emph{Topology Appl.} 40 (1991), no. 2, 189--194. Updated version available at \url{http://www.math.cornell.edu/~hatcher/Papers/TriangSurf.pdf}.

\bibitem[Iv]{Iv} N. V. Ivanov, Mapping class groups, in \emph{Handbook of geometric topology}, 523--633,
North-Holland, 2002.

\bibitem[LS]{LS} R. Lee and R. Szczarba, On the homology and cohomology of congruence subgroups.
 \emph{Invent. Math.}  33  (1976),  no. 1, 15--53.

\bibitem[MSS]{MSS} S. Morita, T. Sakasai, and M. Suzuki, Abelianizations of derivation Lie algebras of free associative algebra and free Lie algebra, preprint (2011), arXiv:1107.3686.

\end{thebibliography}
\end{document}